\documentclass[a4paper,10pt]{article}

\usepackage{amsmath}
\usepackage{amsfonts}
\usepackage{graphicx}
\usepackage{authblk}

\usepackage{bm}

\usepackage[numbers,sort&compress]{natbib} 
\usepackage{xcolor}

\usepackage[
  colorlinks=true,
  linkcolor=blue,
  citecolor=teal,
  urlcolor=teal
]{hyperref}

\advance\oddsidemargin -0.5in
\advance\topmargin -0.5in
\advance\textwidth 1.0in
\advance\textheight 1.0in

\title{A variational approach to estimating the state of a magma reservoir from observed displacement}

\author[1,2]{Shungo Kun Tonoyama}
\author[1]{Atsushi Suzuki\thanks{Corresponding author: {\em atsushi.suzuki.aj@a.riken.jp} }}
\author[1,2]{Takemasa Miyoshi}
\affil[1]{RIKEN Center for Interdisciplinary Theoretical and Mathematical Sciences (iTHEMS) Kobe, Hyogo, 650-0047, Japan}
\affil[2]{RIKEN Center for Computational Science (R-CCS) Kobe, Hyogo, 650-0047, Japan}

\date{}
\begin{document}

\maketitle
\vspace{-2.0em}

\begin{abstract}
  \noindent
We propose a numerical procedure to solve an inverse problem that estimates the state of a magma reservoir from observed surface displacement of a volcano.
  Our variational approach aims to find the minimizer of a cost function consisting of a norm concerning both data and derivative, which evaluates the misfit between the estimated and observed displacement.
  The extremal of the cost function leads to a linear system, to find the stress distribution on the reservoir surface, has very high condition number, but it is feasible to get appropriate solution by using high-precision arithmetic without applying common regularization technique.

\end{abstract}

\section{Introduction}
\noindent
Inverse modeling is widely used in geophysics to infer subsurface processes from surface observations, including the state of magma reservoirs constrained by volcanic deformation. The reliability of such inversions strongly depends on the available observations. Interferometric Synthetic Aperture Radar (InSAR) analysis now provides spatially dense deformation measurements, complementing sparse point observations from Global Navigation Satellite System (GNSS) networks \cite{Massonnet1998}.

In volcanology, analytical source models such as the Mogi model for a spherical pressure source in a homogeneous elastic half-space remain widely used because of their simplicity \cite{Mogi1958}. However, these solutions rely on {strong} assumptions (e.g., a flat free surface) and can be biased when geometric complexity is non-negligible \cite{Cayol1998Effects}.

Here, we propose a variational approach with the adjoint method to estimate the state of a pressure source from surface displacement observations within a linear elastic setting. We formulate the inverse problem as the recovery of the traction (stress) distribution on the reservoir boundary from displacement measurements on the ground surface.
This kind of inverse problem is not treated numerically due to supposed singularity in mathematical point of view. However, an iterative method with high-precision arithmetic can find an appropriate solution, though the linear system derived from our approach still has high condition number.

\section{Model setting}\label{sec:model}

\subsection{Linear elastic deformation model}
\noindent
To generate synthetic surface-displacement data driven solely by a traction prescribed on the magma-reservoir boundary, we implement a three-dimensional linear elastic deformation model by finite element method.
Although viscoelastic formulations have been investigated for volcanic deformation, a linear elastic approximation is appropriate when the deformation is interpreted over sufficiently short time scales \cite{Heap2020}.

\subsection{Geometry and governing equations} 
\noindent
Let $\Omega\subset\mathbb{R}^3$ be a bounded domain representing the host rock. We introduce an interior cavity $\omega\subset\Omega$ whose
boundary $\partial\omega$ represents the magma-reservoir wall. The elastic medium occupies the perforated domain $D := \Omega\setminus\overline{\omega}$ with Lipschitz boundary $\partial{}D$, whose outward unit normal is denoted as $n$.
The exterior boundary $\partial\Omega$ is decomposed into a union of three disjoint parts,
\begin{equation}\label{eq:extBoundary}
\partial\Omega = \Gamma_B \cup \Gamma_W \cup \Gamma_S,
\end{equation}
where $\Gamma_B$ denotes a clamped boundary, $\Gamma_W$ a traction-free boundary, and $\Gamma_S$ the observation (ground) surface. 
\par
The constitutive relation is given by Hooke's law for a homogeneous, isotropic, linear elastic solid,
\begin{equation}\label{eq:hooke}
\bm{\sigma}(u)=\lambda\,\mathrm{tr}\!\left(\bm{\epsilon}(u)\right)I + 2\mu\,\bm{\epsilon}(u),
\end{equation}
where $\lambda$ and $\mu$ are the Lam\'e constants, $I$ is the $3\times 3$ identity tensor, and the infinitesimal strain is
$\bm{\epsilon}(u) = \left(\nabla u + \nabla u^{\mathsf T}\right)/2$\,.
\par
Under the quasi-static assumption, and ignoring the gravity force, the displacement field $u:D\to\mathbb{R}^3$ satisfies kinetic equilibrium,
\begin{equation*} 
 -\nabla\cdot\bm{\sigma}(u) = 0 \qquad \text{in } D\,.
\end{equation*}
with boundary conditions, $u=0$ on $\Gamma_{B}$ and
$\bm{\sigma}(u)n = 0$ on $\Gamma_W$.
\par
In this model, displacement is observed as $u^{\ast}$ on $\Gamma_{S}$ and
traction will vanish, i.e., $\bm{\sigma}(u)n=0$ on the same ground surface boundary. The model finds appropriate traction $\bm{\sigma}(u)n=g$ on the reservoir surface $\partial\omega$.

\begin{figure}[t]
\centering
\includegraphics[scale=0.5]{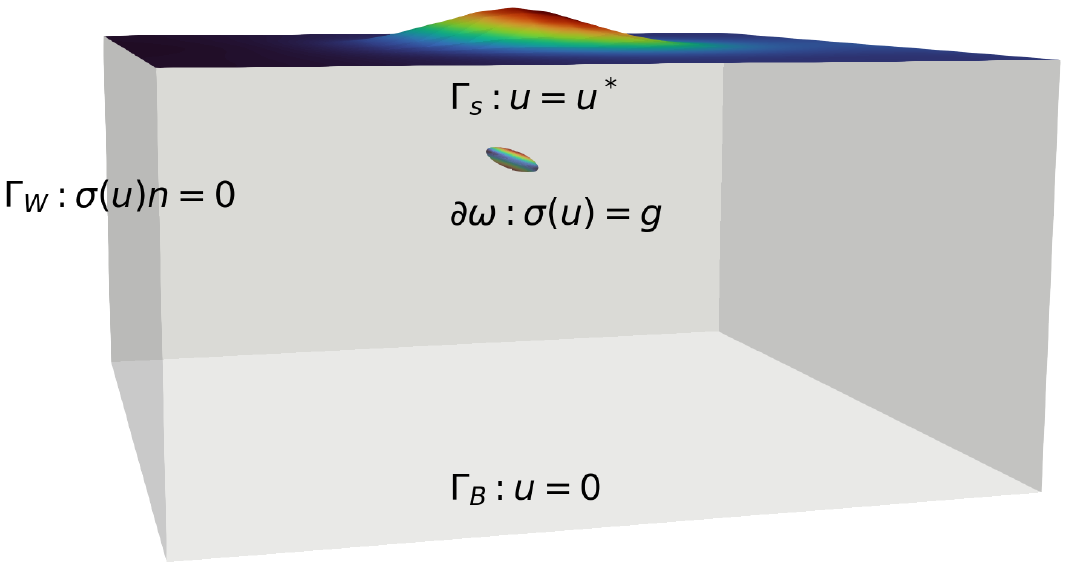}
\caption{Computational domain used in the synthetic experiment.}
\label{Fig:domain}
\end{figure}

\section{Inverse problem formulation}\label{sec:inverse}
\noindent
Since two kinds of boundary conditions on the ground surface boundary, displacement and traction, cannot be set simultaneously, we need to find a way to
estimate traction $g$ on $\partial\omega$ by assuming either $u=u^{\ast}$ or ${\bm\sigma}(u)n=0$ on $\Gamma_{S}$.
The bilinear form is defined as
\begin{equation}\label{eq:bilinear}
a(u,v)=\int_{D}\bm{\sigma}(u):\bm{\epsilon}(v)
\end{equation}
for $u$ and $v$ belonging to the Sobolev space, $H^{1}(D)^{3}$. Precise setting of the solution space will be described in following two subsections.
\subsection{Targeting null traction on the ground surface}
The first method is to give $u=u^{\ast}$ on $\Gamma_{S}$ and $\bm{\sigma}(u)n=g$ on $\partial\omega$ with targeting $\bm{\sigma}(u)n=0$ and $\Gamma_{S}$.
\par
First, we introduce an affine space where Dirichlet data $u=u^{\ast}$ are prescribed on $\Gamma_{S}$,
\begin{equation}\label{eq:affine-space}
  W(u^{\ast}):=\{u\in H^{1}(D)^{3}\,;\, u=0\text{ on }\Gamma_{B}, u=u^{\ast}\text{ on }\Gamma_{S}\}\,.
\end{equation}
The kinetic equilibrium is obtained as a solution of a weak formulation for given stress ${\bm\sigma}(u)n=g$ on $\partial\omega$, to find
$u\in W(u^{\ast})$
\begin{equation*} 
  0=a(u, v)-\langle g, v\rangle_{\partial\omega} \quad \forall v\in W(0)\,.
\end{equation*}
Here $\langle g, v\rangle_{\partial\omega}=\int_{\partial\omega}g\cdot v$ denotes the duality pair of $H^{-1/2}(\partial\omega)^{3}$ and $H^{1/2}(\partial\omega)^{3}$.
\par
We prepare another function space as a superset of $W(u^{\ast})$, where the data on $\Gamma_{S}$ are not fixed,
\begin{equation}\label{eq:test-func-space}
  V:=\{u\in H^{1}(D)^{3}\,;\, u=0\text{ on }\Gamma_{B}\}\,.
\end{equation}
The targeting condition on null traction is written as
\begin{equation*}
  0=\int_{\Gamma_{S}}\bm{\sigma}(u)n\cdot v = \langle \bm{\sigma}(u)n, v\rangle_{\Gamma_{S}}\quad \forall v\in H^{1/2}(\Gamma_{S})^{3}\,.
\end{equation*}
Using integration by parts, we rewrite this surface integration by the domain integration as
\begin{equation*}
\int_{\Gamma_{S}}\bm{\sigma}(u)n \cdot v = \int_{D}\bm{\sigma}(u):\bm{\epsilon}(v)-
  \int_{\partial\omega}g\cdot v \quad u, v\in V
\end{equation*}
and it leads to a problem to find $g\in H^{-1/2}(\partial\omega)$ and $\hat{u}(g)\in W(u^{\ast})$ satisfying
\begin{equation}\label{eq:underdet}
  0=a(\hat{u}(g), v)-\langle g, v\rangle_{\partial\omega} \quad \forall v\in V\,.
\end{equation}
Forthcoming subsection \ref{seq:matrix-normal} shows matrix representation of this weak formulation and a way to solve the linear system via a normal equation in algebraic manner.
\subsection{Minimizing misfit on the ground surface}
We would like to introduce the second method
  by assuming $\bm{\sigma}(u)n=0$ and to target observed data $u^{\ast}$ on
  $\Gamma_{S}$, which derives a direct weak formulation on the reservoir wall.
Let $u(g)\in V$ satisfying the kinetic equilibrium,
\begin{equation}\label{eq:kinetic-state}
a(u(g), v)=\langle g, v\rangle\quad \forall v \in V\,.
\end{equation}
We define a cost function to minimize the data-misfit,
\begin{equation}\label{eq:cost}
J(g)=\frac12\|u(g)-u^\ast\|_{H^{1/2}(\Gamma_S)^{3}}^2.
\end{equation}
Here, the $H^{1/2}(\Gamma_S)^{3}$-norm is defined via a harmonic extension of data $w$ on $\Gamma_{S}$ to $\widetilde w$ in $D$, by solving a problem to find $\widetilde w\in W(w)$ satisfying
\begin{equation}\label{eq:harmonic-ext}
 a({\widetilde w}, v)=0 \quad \forall v \in W(0)
\end{equation}
and $\|w\|_{H^{1/2}(\Gamma_{S})^{3}}^{2}=a(\widetilde w, \widetilde w)$.
By putting another extended function $\widetilde u$ of $u$ into \eqref{eq:harmonic-ext}, we have
\begin{equation*}
a(\widetilde{w}, \widetilde{u})=\langle \bm{\sigma}(\widetilde{w})n, u\rangle_{\Gamma_{S}}
\end{equation*}
since $u$ does not vanish on $\Gamma_{S}$. This relation guarantees the cost function is expressed by a duality pair of the data misfit and whose traction.
\par
Now we formulate a problem to find traction $g\in H^{-1/2}(\partial\omega)^{3}$ that minimizes \eqref{eq:cost} with the constraint that is given as \eqref{eq:kinetic-state}.
To compute the gradient of $J(g)$, we introduce a Lagrangian
for $u, v\in V$ in addition to $g$,
\begin{equation}\label{eq:Lagrangian}
  \mathcal{L}(u,v;g)=J(g)- a(u,v) + \langle g,v\rangle_{\partial\omega}\,.
\end{equation}
We note that, by using the harmonic extensions of $u$ and $u^{\ast}$ as
$\widetilde{u}$ and $\widetilde{u^{\ast}}$, respectively, $J(g)$ is expressed as
\begin{equation*}
  J(g)=\frac12{}a(\widetilde u - \widetilde{u^{\ast}},
  \widetilde u - \widetilde{u^{\ast}})
  =\frac12{}\langle\bm{\sigma}(\widetilde{u(g)}-\widetilde{u^{\ast}})n, u(g)-u^{\ast}\rangle_{\Gamma_{S}}\,.
\end{equation*}
The standard adjoint method\cite{Azegami2020} by solving two systems as the state problem for $u(g)$ by fixing $g$ and the adjoint problem $v(g)$ by fixing $g$ and $u(g)$ helps to calculate  gradient of $J(g)=\mathcal{L}(u(g), v(g); g)$ by finding the saddle point of the Lagrangian.
A variation of \eqref{eq:Lagrangian} by $v+\delta v$ leads to 
the state problem, which is nothing but the kinetic equilibrium \eqref{eq:kinetic-state}\,. 
And then, another variation of \eqref{eq:Lagrangian} by $u(g)+\delta u$ with $\delta u\in V$ leads to the adjoint problem, which is stated as to find $v\in V$ satisfying
\begin{equation}\label{eq:adjoint}
  a(\delta u, v)=\langle \delta u, \bm{\sigma}(\widetilde{u(g)})n - \bm{\sigma}(\widetilde{u^{\ast}})n \rangle_{\Gamma_{S}} \quad \forall \delta u\in V\,.
\end{equation}
By estimating variation of the state and adjoint solutions $u(g)$ and $v(g)$ with respect to $g+\delta g$, we have 
\begin{equation*} 
\delta J(g)[\delta g]=\langle \delta g,\ v(g)\rangle_{\partial\omega}.
\end{equation*}
An extremal of $J(g)$ is attained by a following variational problem to find $g\in H^{-1/2}(\partial\omega)^{3}$ satisfying
\begin{equation}\label{eq:variational}
  \langle \delta g,\ \bar v(g)\rangle_{\partial\omega}=
  \langle \delta g,\ v^{\ast}\rangle_{\partial\omega}\ \ \forall \delta g\in H^{-1/2}(\partial\omega)^{3}\,.
\end{equation}
Here $\bar v(g)$ is the solution of \eqref{eq:adjoint} with the harmonic extension $\widetilde{u(g)}$ from $u(g)$ and $v^{\ast}$ calculated from ${u^{\ast}}$, respectively. 
The quadratic form of $J(g)$ and $v(g+\delta g)=v(g)+{\bar v}(\delta g)$ leads to
 $J(g+\delta g)-J(g)=(1/2)\langle \delta g, {\bar v}(\delta g)\rangle_{\partial\omega}$, which shows positivity guarantees the unique minimizer.

\section{Matrix representation}\label{sec:algo}
\noindent
In this section we present the matrix formulation by finite element discretization of the weak form of the kinetic equilibrium, which will be applied to two methods. 
\subsection{Block decomposition of the stiffness system}
Let $K$ be the stiffness matrix associated with the linear elasticity problem on
$D=\Omega\setminus\overline{\omega}$.
The discrete system for the kinetic equilibrium \eqref{eq:kinetic-state} reads
\begin{equation}\label{eq:blockK}
\begin{bmatrix}
K_{11} & K_{12} & 0\\
K_{21} & K_{22} & K_{23}\\
0      & K_{32} & K_{33}
\end{bmatrix}
\begin{bmatrix}
\mathbf{u}_{1}\\ \mathbf{u}_{2}\\ \mathbf{u}_{3}
\end{bmatrix}
=
\begin{bmatrix}
\mathbf{g}_{1}\\ 0\\ 0
\end{bmatrix}\,.
\end{equation}
Here, degrees of freedom is ordered as $\Lambda=\Lambda_{1}\oplus\Lambda_{2}\oplus\Lambda_{3}\oplus\Lambda_{B}$, corresponding to 1) nodes on $\partial\omega$, 2) internal nodes in $D\cup\Gamma_{W}$, 3) nodes on $\Gamma_{S}$, and B) nodes on $\Gamma_{B}$.
The vector $\mathbf{u}=(\mathbf{u}_{1},\mathbf{u}_{2},\mathbf{u}_{3})$ follows this order, and
$\mathbf{g}_{1}$ is the discrete 
traction load vector induced by the unknown traction $g$ on $\partial\omega$ with surface integration. 
Let $\mathbf{u}_{3}^\ast$ denote the observed surface displacement on $\Gamma_S$.

\subsection{Surface response map and Schur complements}\label{seq:matrix-normal}
The weak problem \eqref{eq:underdet} is expressed by the linear system \eqref{eq:blockK} with setting $\mathbf{u}_{3}=\mathbf{u}_{3}^\ast$.
However, unknown vectors $\mathbf{u}_{1}$, $\mathbf{u}_{2}$ and $\mathbf{g}_{1}$ are located in both sides of the system.
Eliminating $\mathbf{u}_{1}$ and $\mathbf{u}_{2}$ yields the linear mapping from $\mathbf{g}_{1}$ to the surface displacement $\mathbf{u}_{3}$,
\begin{equation}\label{eq:Amap}
\mathbf{u}_{3} = A\,\mathbf{g}_{1}, \qquad A := R_{3}K^{-1}R_{1}^{\mathsf T}.
\end{equation}
Using block elimination, we introduce the Schur complements
$S_{22}=K_{22}-K_{21}K_{11}^{-1}K_{12}$ and
$S_{33}=K_{33}-K_{32}S_{22}^{-1}K_{23}$,
and we define
\begin{equation}\label{eq:B31}
B_{31}:=K_{32}S_{22}^{-1}K_{21}K_{11}^{-1}.
\end{equation}
Then the surface response \eqref{eq:Amap} is written as
\begin{equation}\label{eq:u3}
\mathbf{u}_{3} = S_{33}^{-1}B_{31}\,\mathbf{g}_{1},
\end{equation}
which will be under-determined when $\#\Lambda_{1} >\#\Lambda_{3}$.
The normal equation associated with matrix $A=S_{33}^{-1}B_{31}$,
\begin{equation}\label{eq:normal}
A^{\mathsf T}A\,\mathbf{g}_{1} = A^{\mathsf T}\mathbf{u}_{3}^\ast
\end{equation}
can be solved by conjugate gradient or GMRES method in $\text{Im}A^{\mathsf T}$,
but the normal equation will scale the condition number as double.

\subsection{Variational equation via harmonic extension}
Here, we present a matrix expression of the variational problem \eqref{eq:variational} on discrete traction data $\mathbf{g}_{1}$.
  Since the $H^{1/2}(\Gamma_S)^{3}$-inner product is defined via a harmonic-extension,
the solution of the adjoint problem \eqref{eq:adjoint} in discrete sense, $[\mathbf{u}_{1}, \mathbf{u}_{2}]$ for prescribed $\mathbf{u_{3}}$ is given by
\begin{equation}\label{eq:hext}
\begin{bmatrix}
K_{11} & K_{12}\\
K_{21} & K_{22}
\end{bmatrix}
\begin{bmatrix}
\mathbf{u}_{1}\\ \mathbf{u}_{2}
\end{bmatrix}
=
\begin{bmatrix}
0\\ -K_{23}\mathbf{u}_{3}
\end{bmatrix}.
\end{equation}
This leads to a relation from the surface to the interface,
\begin{equation} \label{eq:B13}
\mathbf{u}_{1} = B_{13}\,\mathbf{u}_{3},\qquad
B_{13}:=K_{11}^{-1}K_{12}S_{22}^{-1}K_{23}.
\end{equation}
Combining with \eqref{eq:u3}, we have the matrix representation,
\begin{equation}\label{eq:variational-discrete}
B_{13}S_{33}^{-1}B_{31}\,\mathbf{g}_{1} = B_{13}\mathbf{u}_{3}^\ast.
\end{equation}
Here, the right-hand side corresponds to the map from the surface data $u^{\ast}$ to the trace of its harmonic extension $v^{\ast}$.
The map from $g$ to $\bar v(g)$ is represented by 
$B_{13}S_{33}^{-1}B_{31}$ in the left-hand side of \eqref{eq:variational-discrete}.
Since the Schur complement $S_{33}$ is invertible, $B_{13}S_{33}^{-1}B_{31}$ is symmetric positive semi-definite, which corresponds to singularity of \eqref{eq:variational}.
\par
In our implementation, we primarily work with \eqref{eq:variational-discrete}, which is better in the matrix property than \eqref{eq:normal}, whose coefficient matrix $A^{\mathsf T}A$ is expressed as $B_{13}S_{33}^{-1}S_{33}^{-1}B_{31}$.
\par
For the case when more degrees of freedom on $\partial\omega$ 
accommodate than ones on $\Gamma_{S}$, 
we need to keep in mind $\#\Lambda_{1} >\#\Lambda_{3}$, which leads to a semi-definite system even in algebraic sense. By aggregating neighboring finite element nodes and setting the same data on paired nodes, such pseudo singularity is eliminated.

\section{Regularization technique}
The Tikhonov regularization\cite{Gockenbach2016} is widely used to solve least squares minimization problem with rectangular matrix.
The under-determined system \eqref{eq:Amap} is reformulated as a minimization problem to find $\mathbf{g}_{1}$ attaining $\min_{\mathbf{g}_{1}}||A\,\mathbf{g}_{1}-\mathbf{u}_{3}^{\ast}||^{2}$.
By preparing a matrix $L$ satisfying $\text{Ker}A\cap\text{Ker}L=\emptyset$ and a positive parameter $\varepsilon >0$, perturbed function to evaluate the misfit is expressed as
\begin{equation*}
  ||A\,\mathbf{g}_{1}-\mathbf{u}_{3}^{\ast}||^{2}+\varepsilon||L\,\mathbf{g}_{1}||^{2}\,.
\end{equation*}
\par
Let us apply this technique to our adjoint-based method. 
We set a new cost function by adding kinetic energy of the state solution $u(g)$ with parameter $\varepsilon>0$,
\begin{equation}\label{eq:cost-tikhonov}
  \widetilde{J(g)}=J(g)+\varepsilon\langle g, u(g)\rangle_{\partial\omega}\,.
\end{equation}
The adjoint problem is formulated by adding $\varepsilon\langle g, \delta u\rangle$ to the right hand side of \eqref{eq:adjoint}, which leads to a variational problem to find $g$ with using the same $\bar v(g)$ in \eqref{eq:variational} as
\begin{equation}\label{eq:variational-tikhonov}
  \langle \delta g,\ {\bar v}(g)\rangle_{\partial\omega}
+2\varepsilon\langle \delta g, u(g)\rangle_{\partial\omega}=  
  \langle \delta g,\ v^{\ast}\rangle_{\partial\omega}\,.
\end{equation}
Here, $g$ and $\delta g$ are taken in 
$H^{-1/2}(\partial\omega)^{3}$ as same as \eqref{eq:variational}\,.
A discretized form using the same symbols of \eqref{eq:variational-discrete} reads
\begin{equation}\label{eq:Tikhonov-discrete}
(B_{13}S_{33}^{-1}B_{31}+2\varepsilon R_{1}K^{-1}R_{1}^{T})\,\mathbf{g}_{1} = B_{13}\mathbf{u}_{3}^\ast.
\end{equation}
Though \eqref{eq:variational-tikhonov} has a unique solution, which attains the minimizer of \eqref{eq:cost-tikhonov}, it strongly depends on $\varepsilon$.
\par
We will compare solutions of \eqref{eq:variational-discrete} and \eqref{eq:Tikhonov-discrete} and show an advantage of solution of \eqref{eq:variational-discrete} by iterative solver in the image space with high-precision arithmetic.

\begin{figure}[t]
\centering
\includegraphics[scale=0.55]{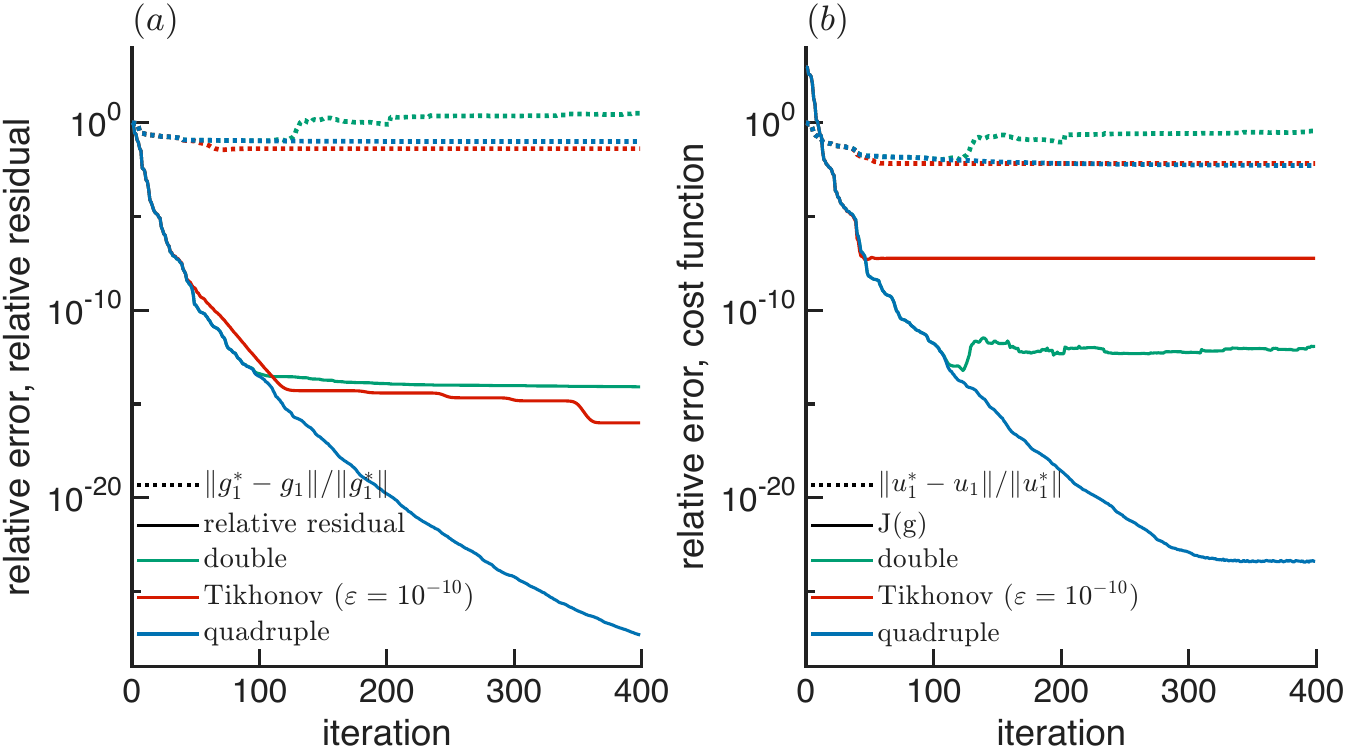}
\caption{
Convergence history of the adjoint-based inversion. ($a$) Relative residual norm and traction reconstruction error on $\partial\omega$. ($b$) Cost function $J(g)$ and displacement error on $\Gamma_s$. Green, red, and blue curves denote history calculated by double-precision, Tikhonov regularization, and quadruple-precision, respectively.
}
\label{Fig:itr}
\end{figure}

\section{Result and discussion}\label{sec:result}
\noindent
We validated the proposed adjoint-based inversion using a synthetic benchmark in a three-dimensional elastic domain with a Gaussian-smoothed free surface and an ellipsoidal cavity $\omega$ representing a magma reservoir (Fig.~\ref{Fig:domain}). A spatially varying traction $g^{*}$ prescribed on the 
boundary $\partial\omega$ was 
set to follow a Gaussian distribution, and
by taking $\lambda=\mu=1$, the resulting displacement field $u^{*}$ on the 
surface $\Gamma_{S}$
was treated as the observation. 
\par
For finite element computation, P1 element is used and the total degrees of freedom is 226,953 and 34,029 / 29,136 DOFs on $\partial\omega$ and $\Gamma_{S}$, respectively.
FreeFEM is used to generate the stiffness matrix from given unstructured mesh created by Gmsh.
GMRES iteration with double/quadruple precision is used to solve \eqref{eq:variational-discrete}.

\begin{figure}[t]
\centering
\includegraphics[scale=0.35]{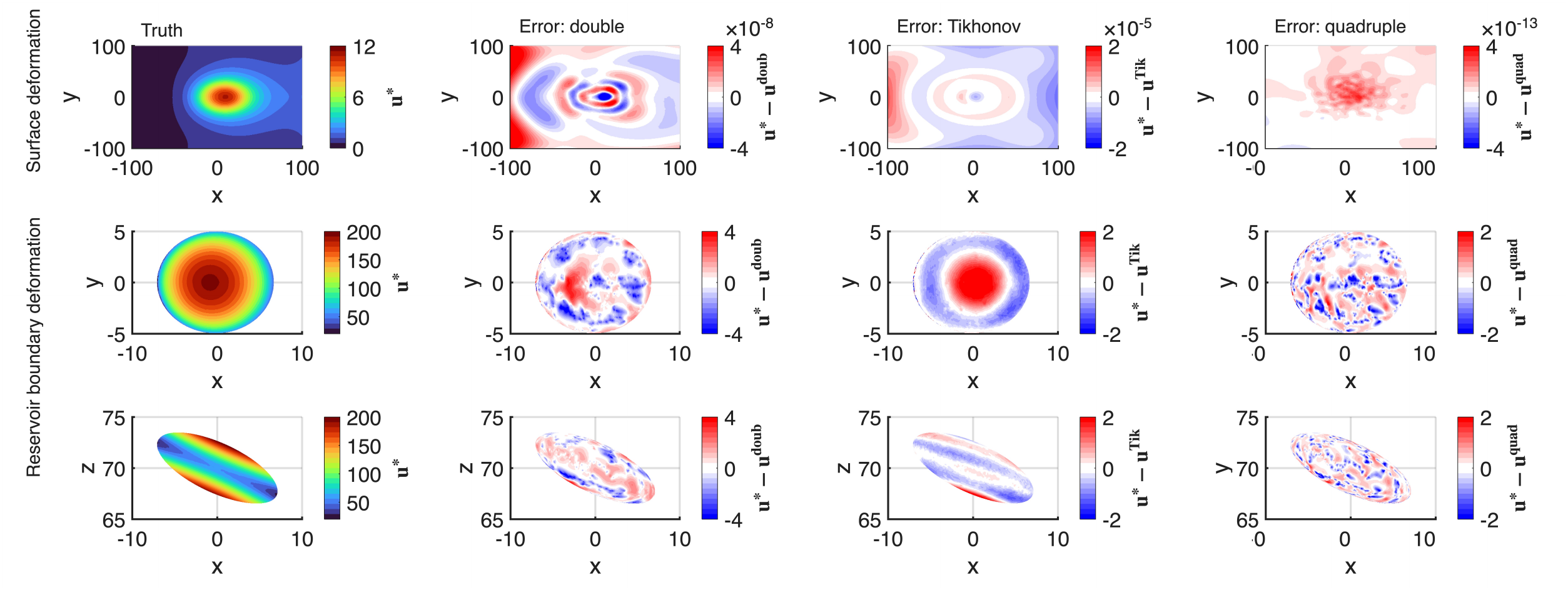} 
\caption{
True displacement and reconstruction errors on the ground surface $\Gamma_S$ and the reservoir boundary $\partial\omega$. Columns show, from left to right, the true displacement $\mathbf{u}^{*}$, the errors by double-precision $\mathbf{u}^{*}-\mathbf{u}^{\mathrm{doub}}$ after 100 iterations, by Tikhonov-regularization $\mathbf{u}^{*}-\mathbf{u}^{\mathrm{Tik}}$ after 100 iterations, and by quadruple-precision $\mathbf{u}^{*}-\mathbf{u}^{\mathrm{quad}}$ after 300 iterations. The iteration numbers were selected based on the stability in Fig. 2. The top row shows displacement on $\Gamma_S$; the middle and bottom rows show bottom and side views on $\partial\omega$.
}
\label{Fig:err}
\end{figure}

We solved the state
problem up to 400 iterations (Fig.~\ref{Fig:itr}). In double-precision arithmetic, the residual decreases rapidly at early iterations but begins to saturate after approximately $\sim$100 iterations, and the error 
becomes even larger beyond 150 iterations (Fig.~\ref{Fig:itr}~($a$)). 
In contrast, quadruple-precision arithmetic, both in the global GMRES iteration and in a direct solver~\cite{Suzuki2022} for the kinetic equilibrium system, enables the computation to proceed further.
Iteration by the Tikhonov regularization with $\varepsilon=10^{-10}$ takes different history path after 50 iterations but slightly better residual than ``quadruple'' is obtained (Fig.~\ref{Fig:itr}~($a$)).
The cost function $J(g)$ decreases up to 300 iterations by ``quadruple'', but it stagnates at around $10^{-8}$ by ``Tikhonov''
(Fig.~\ref{Fig:itr}~($b$)).
The reconstructed displacement on $\partial\omega$ decreases at the same level in ``quadruple'' and ``Tikhonov''.
Spatial comparisons further confirm the effectiveness of the inversion. The reconstructed traction reproduces the imposed Gaussian-like pattern on both $\Gamma_{S}$ and $\partial\omega$, and the remaining residual shows no coherent structure in ``quadruple'' (Fig.~\ref{Fig:err}).
On the contrary, substantial perturbation is observed in ``Tikhonov''
on the bottom side of the magma reservoir $\partial\omega$ and error of the reconstructed displacement on the surface $\Gamma_{S}$ is not completely negligible.
\par
In this experiment, the final relative error on the reservoir boundary is approximately $\sim$1\%, demonstrating that the adjoint-based framework can recover physically meaningful boundary loading from surface displacement.
\par
An improvement by combination of the regularization with sufficiently small perturbation and high-precision arithmetic will be investigated in future work.
\par
These results suggest that the proposed methodology will provide a useful framework for inferring stress concentrations on magma reservoirs from surface displacement, where future eruptive activity could be caused.

\section{Acknowledgements}
\noindent
This work was supported by the funding of the RIKEN TRIP initiative (Prediction Science); Fugaku Kodoka (ra000007); JAXA EORA4 (ER4A2N524, ER4MAF004); the COE research grant in computational science from Hyogo Prefecture and Kobe City through Foundation for Computational Science; JSPS KAKENHI (JP24H00021); JST (JPMJSA2109, JPMJCR24Q3); and the UK ARIA (FPCW-PR01-P007).

\bibliographystyle{unsrtnat}
\bibliography{reference}

\end{document}